\documentclass[]{article}
\usepackage{geometry}                
\usepackage[sc]{mathpazo}
\usepackage{graphicx}
\usepackage{enumerate}
\usepackage{amsmath,amssymb,amsthm}
\usepackage{xcolor}
\usepackage{framed}
\usepackage{mdframed}
\usepackage{soul}
\usepackage[colorlinks,
      plainpages=false,
      pdfpagelabels=true,
     	pdfpagemode=UseOutlines,
      pdfpagelayout=SinglePage,
			bookmarks=false,
			colorlinks=true,
      hyperfootnotes=false,
			linkcolor=blue!50!black,
			citecolor=green!40!black]{hyperref}

\newtheorem{remark}{Remark}
\newtheorem{lemma}[remark]{Lemma}

\newtheorem{theorem}[remark]{Theorem}

\newcommand{\hlc}[2][yellow]{ {\sethlcolor{#1} \hl{#2}} }

\newenvironment{myproof}[1][Proof]{\begin{shaded}\vspace*{-1em}\begin{proof}[#1]}{\end{proof}\vspace*{-1em}\end{shaded}}

\colorlet{shadecolor}{gray!10}

\newcommand{\R}{\mathbb{R}}

\newcommand{\mat}[1]{\mathbf{#1}}
\newcommand{\set}[1]{\mathcal{#1}}

\title{Split optimal policy iteration for LQR problems}
\author{P\'eter Koltai}
\date{\today}                                           

\begin{document}
\maketitle

\begin{abstract}
This technical report is concerned with the convergence properties of what we call the \emph{split optimal policy iteration} for coupled LQR problems; see section~\ref{sec:twosubsys} below. Interestingly, the iteration shows different convergence behavior for continuous and discrete time systems: while global convergence holds for both cases, we have local quadratic convergence for the continuous time case, but only linear convergence for the discrete time case---even though quadratic convergence is retained in the limit as the coupling between the subsystems vanishes.
\end{abstract}

\tableofcontents

\section{Introduction}

We consider an iterative solution method for optimal control problems for time-invariant coupled systems. For time-invariant (autonomous) systems the optimal controller, if it exists, is given by a \emph{feedback} (also called \emph{policy}), i.e.\ a mapping from the state into the admissible control space. We consider the solution to the optimal control problem being equal to finding this policy. The system structure considered here is such that each subsystem is characterized by its own state and control variables, and the evolution of the state variables is coupled through the dynamics. For such systems one could attempt to solve the optimal control problem by optimizing the policy of just one subsystem while keeping the policies of the other subsystems fixed, then doing the same for the next subsystem, and so on. Since in each step the optimal policy with respect to the current policies of the other subsystems is chosen by the current subsystem, we call this split optimal policy iteration. The process is exemplified for LQR problems with two subsystems in section~\ref{sec:twosubsys}.

\smallskip
The purpose of this technical report is the \emph{qualitative and quantitative} analysis of the split optimal policy iteration applied to continuous and discrete time LQR problems. The motivation is rather phenomenological, since it is not immediate to the author whether this procedure has any numerical advantage over solving the algebraic Riccati equation of the full problem directly. However, if the here observed \emph{fast} convergence for weakly coupled systems carries over to nonlinear systems as well, then this could open pathways for the design of centralized (sub-)optimal controllers for general nonlinear systems, as discussed in~\cite{KoJu14b}.

\smallskip
The idea of split iterative optimization is by no means new. For instance, iterative optimization over reduced variables can be considered for minimizing scalar functions $f:\R^n\to\R$ whenever the structure of $f$ is such that it is easier to solve $\min_{x_i} f(x_1,\ldots,x_n)$ for any $i=1,\ldots,n$ than the full problem, and the cyclical repetition of this procedure is expected to converge to a reasonable solution~\cite{BeHa02}. Note the difference to the setting considered here: in optimal control one optimizes the optimal value function in every state of the state space, hence we are dealing with multiple objective functions. However, Bellman's optimality principle~\cite{Bell57} states that these objectives are not concurrent; in order to be optimal, every end-part of a path has to be optimal itself. This would suggest that instruments from this field are likely to be able to be carried over to our setting. Indeed, the monotonicity in Lemma~\ref{lem:monotony_cLQR} can be seen as an example for this.

\smallskip
Another property of the iterative procedure described here is that for general nonlinear systems it is not guaranteed to converge to the optimal solution---just as the iteration applied to minimize $f:\R^n\to\R$ can converge to a local but not global minimum. However, a fixed point of the iteration has a game-theoretical interpretation: it is a \emph{Nash equilibrium}~\cite{OsRu94}. In a Nash equilibrium the policy of each player (here: the subsystem controller) is optimal with respect to the policies of the other players.\footnote{The reader might be irritated by the fact that in our coupled optimal control setting the subsystems follow a joint \emph{cooperative} goal, while game theory usually considers \emph{competitive} situations. However, the ``definition'' of a Nash equilibrium given here does not include the players' intentions, it only assumes they are optimal with respect to each other. This covers both situations.} In game-theoretical terms our procedure is a \emph{best-response strategy} iteration, or a kind of \emph{fictitious play}, which has been proposed to compute a Nash equilibrium of a game~\cite{Bro51,Rob51,MoSh96a,MoSh96b}. It has to be noted, that games often only possess \emph{mixed} Nash equilibria, which are randomized strategies. Optimal control problems (even with stochastic components, as \emph{Markov decision processes}), however, have deterministic optimal policies. But coupled optimal control problems can be rewritten as static (collaborative) games, see e.g.\ the lecture notes~\cite{Johari}. It is not clear to the author whether, and to which extent, does this connection allow to carry over the convergence results from game theory to the present setting.

\smallskip
As the aim of the subsystems is to optimally accomplish a joint goal, it is not surprising that very similar ideas occur in \emph{multiagent systems} and \emph{reinforcement learning} too~\cite{Litt01,Meuleau98solvingvery,Lauer00analgorithm,Guestrin01multiagentplanning}. Certainly, the idea of separability can be found in the \emph{dynamic programming} literature too~\cite{Ber07}.

\smallskip
To facilitate the navigation in this report, the \fbox{main theorems are framed} and\hlc[shadecolor]{all proofs have a gray background} (and can be skipped if one is only interested in the main results).

\vfill

\section{LQR problems: notation}
\subsection{Continuous time}

Linear-quadratic regulator (LQR) problems arise in the context when the linear system
\begin{equation}
\dot x(t) = \mat{A}x(t) + \mat{B}u(t)
\label{eq:linsysC}
\end{equation}
shall be controlled to the origin in an optimal way. Here, $\mat{A}\in\R^{m\times m}$ denotes the system matrix and $\mat{B}\in\R^{m\times r}$ denotes the input matrix. Further, let $\mat{Q}\in\R^{m\times m}$ and $\mat{R}\in\R^{r\times r}$ be symmetric positive definite matrices. The control task is to find a function $u:[0,\infty)\to \set{U}$, generating a trajectory~$\{x(t)\}_{t\ge 0}$ through~\eqref{eq:linsysC}, such that the accumulated costs
\[
\int_0^{\infty} x(t)^{\rm T}\mat{Q}x(t) + u(t)^{\rm T}\mat{R}u(t)\,dt
\]
are minimal. By construction, this already implies that ${x(t)\to 0}$ as~$t\to\infty$. It turns out that the optimality principle is equivalent to the \textit{continuous algebraic Riccati equation}
\begin{equation}
\mat{A}^{\rm T}\mat{P} + \mat{P}\mat{A} - \mat{P}\mat{B}\mat{R}^{-1}\mat{B}^{\rm T}\mat{P} + \mat{Q} = \mat{0}
\label{eq:CARE}
\end{equation}
where the unique symmetric positive definite solution~$\mat{P}^{\rm opt}$ of this equation yields the optimal value function ${V(x) = x^{\rm T}\mat{P}^{\rm opt}x}$. Moreover, the optimal feedback is given by
\begin{equation}
\mu(x) = \mat{F}^{\rm opt}x = -\mat{R}^{-1}\mat{B}^{\rm T}\mat{P}^{\rm opt}x.
\label{eq:contOptF}
\end{equation}

\subsection{Discrete time}

Linear-quadratic regulator problems arise in the context when the linear system
\begin{equation}
x(k+1) = \mat{A}x(k) + \mat{B}u(k)
\label{eq:linsysD}
\end{equation}
shall be controlled to the origin in an optimal way. Here, $\mat{A}\in\R^{m\times m}$ denotes the system matrix and $\mat{B}\in\R^{m\times r}$ denotes the input matrix. Further, let $\mat{Q}\in\R^{m\times m}$ and $\mat{R}\in\R^{r\times r}$ be symmetric positive definite matrices. The control task is to find a sequence $\{u(k)\}_{k\ge0}$, generating a sequence~$\{x(k)\}_{k\ge 0}$ through~\eqref{eq:linsysD}, such that the accumulated costs
\[
\sum_{k\ge 0} x(k)^{\rm T}\mat{Q}x(k) + u(k)^{\rm T}\mat{R}u(k)
\]
are minimal. By construction, this already implies that ${x(k)\to 0}$ as~$k\to\infty$. It turns out that the optimality principle is equivalent to the \textit{discrete algebraic Riccati equation}
\begin{equation}
\mat{P} = \mat{A}^{\rm T}\mat{P}\mat{A} - \mat{A}^{\rm T}\mat{P}\mat{B}\left(\mat{R}+\mat{B}^{\rm T}\mat{P}\mat{B}\right)^{-1}\mat{B}^{\rm T}\mat{P}\mat{A} + \mat{Q},
\label{eq:DARE}
\end{equation}
where the unique symmetric positive definite solution~$\mat{P}^{\rm opt}$ of this equation yields the optimal value function ${V(x) = x^{\rm T}\mat{P}^{\rm opt}x}$. Moreover, the optimal feedback is given by
\begin{equation}
\mu(x) = \mat{F}^{\rm opt}x = -\left(\mat{R}+\mat{B}^{\rm T}\mat{P}\mat{B}\right)^{-1}\mat{B}^{\rm T}\mat{P}^{\rm opt}\mat{A}x.
\label{eq:discrOptF}
\end{equation}

\section{Split optimal policy iteration for LQR problems}
\subsection{Two subsystems}	\label{sec:twosubsys}

In order to get a better intuition about how the split optimal policy iteration works, we show an update step for the case of two subsystems. For this, let us partition the involved matrices $\mat{A}$, $\mat{B}$, $\mat{Q}$, and $\mat{R}$ into blocks according to the subsystem decomposition, i.e.\
\[
\begin{array}{c}
\mat{A} = \begin{pmatrix} \mat{A}_{11} & \mat{A}_{12} \\ \mat{A}_{21} & \mat{A}_{22} \end{pmatrix},\qquad
\mat{B} = \begin{pmatrix} \mat{B}_{11} & \mat{B}_{12} \\ \mat{B}_{21} & \mat{B}_{22} \end{pmatrix},\\[1em]

\mat{Q} = \begin{pmatrix} \mat{Q}_1 & \mat{0} \\ \mat{0} & \mat{Q}_2 \end{pmatrix},\qquad
\mat{R} = \begin{pmatrix} \mat{R}_1 & \mat{0} \\ \mat{0} & \mat{R}_2 \end{pmatrix}.
\end{array}
\]
Assuming that some feedback matrix
\[
\mat{F} = \begin{pmatrix} \mat{F}_{11} & \mat{F}_{12} \\ \mat{F}_{21} & \mat{F}_{22} \end{pmatrix}
\]
is given (partitioned the same way), we now show how the feedback~$\nu_1$ of the first subsystem is updated by the split optimal policy iteration algorithm. The update works analogously for the second subsystem as well. First, note that $\nu_1(x) = \begin{pmatrix} \mat{F}_{11} & \mat{F}_{12}\end{pmatrix} x$ and $\nu_2(x) = \begin{pmatrix} \mat{F}_{21} & \mat{F}_{22}\end{pmatrix} x$. Since~$\nu_2$ is fixed during this update step, we have to merge it into the matrices~$\mat{A}$ and~$\mat{Q}$. We obtain a new system matrix
\[
\mat{A}^{(1)} = \mat{A} + \begin{pmatrix} \mat{B}_{12} \\ \mat{B}_{22} \end{pmatrix}\begin{pmatrix}\mat{F}_{21} & \mat{F}_{22} \end{pmatrix},
\]
and a new state cost matrix
\[
\mat{Q}^{(1)} = \mat{Q} + \begin{pmatrix} \mat{F}_{21} & \mat{F}_{22}\end{pmatrix}^{\rm T}\mat{R}_2 \begin{pmatrix} \mat{F}_{21} & \mat{F}_{22}\end{pmatrix}.
\]
The control- and control cost matrices reduce to
\[
\mat{B}^{(1)} = \begin{pmatrix} \mat{B}_{11} \\ \mat{B}_{21}\end{pmatrix},\ \text{and}\ \mat{R}^{(1)} = \mat{R}_1.
\]
Moreover, we obtain an LQR problem with matrix quadruple $\mat{A}^{(1)}$, $\mat{B}^{(1)}$, $\mat{Q}^{(1)}$, and $\mat{R}^{(1)}$. Solving this LQR problem via \eqref{eq:CARE}  and \eqref{eq:contOptF} or \eqref{eq:DARE}  and \eqref{eq:discrOptF} yields $\nu_1^{\rm new}(x) = \begin{pmatrix} \mat{F}^{\rm new}_{11} & \mat{F}^{\rm new}_{12}\end{pmatrix} x$.

Now, that the new feedback matrix
\[
\mat{F}^{\rm new} =  \begin{pmatrix} \mat{F}^{\rm new}_{11} & \mat{F}^{\rm new}_{12}\\ \mat{F}_{21} & \mat{F}_{22}\end{pmatrix}
\]
is obtained, the iteration continues with the update of the $2^{\rm nd}$ subsystem with~$\mat{F}^{\rm new}$ replacing~$\mat{F}$.

\subsection{Arbitrary number of subsystems}	\label{sec:update_general}

To generalize the approach described in the previous section, consider an LQR problem given by the matrices $\mat{A}$, $\mat{B}$, $\mat{Q}$, and $\mat{R}$, where $\mat{Q}$ and $\mat{R}$ are block diagonal symmetric positive definite. Let $\mat{\Pi}_i$, $i=1,\ldots,n$ with~$n$ being the number of subsystems, be a block column matrix consistent with the decomposition of $\mat{R}$, such that the $i^{\rm th}$ block entry is the identity matrix~$\mat{I}$ and all others are~$\mat{0}$, i.e.\ $\mat{\Pi}_i$ is of the form $(\mat{0},\ldots,\mat{0},\mat{I},\mat{0},\ldots,\mat{0})^{\rm T}$. Analogously as before, this yields
\begin{equation}
\label{eq:matrixredef}
\begin{aligned}
\mat{A}^{(i)} &= \mat{A} + \mat{B} (\mat{I} - \mat{\Pi}_i\mat{\Pi}_i^{\rm T}) \mat{F}, \\
\mat{B}^{(i)} &= \mat{B}\mat{\Pi}_i, \\
\mat{Q}^{(i)} &= \mat{Q} + \mat{F}^{\rm T}(\mat{I} - \mat{\Pi}_i\mat{\Pi}_i^{\rm T})\mat{R}(\mat{I} - \mat{\Pi}_i\mat{\Pi}_i^{\rm T})\mat{F}, \\
\mat{R}^{(i)} &= \mat{\Pi}_i^{\rm T}\mat{R}\mat{\Pi}_i.
\end{aligned}
\end{equation}
Solving the algebraic Riccati equations with these matrices gives $\mat{P}^{\rm new}$, and the feedback matrix is updated by
\begin{eqnarray*}
\mat{\Pi}_i^{\rm T}\mat{F}^{\rm new} & = &
\left\{\begin{array}{ll}
-{\mat{R}^{(i)}}^{-1}{\mat{B}^{(i)}}^{\rm T}\mat{P}^{\rm new}, & \text{continuous time}\\
-\left(\mat{R}^{(i)}+{\mat{B}^{(i)}}^{\rm T}\mat{P}^{\rm new}\mat{B}^{(i)}\right)^{-1}{\mat{B}^{(i)}}^{\rm T}\mat{P}^{\rm new}\mat{A}^{(i)}, & \text{discrete time}
\end{array}
\right\},\\
(\mat{I} - \mat{\Pi}_i^{\rm T})\mat{F}^{\rm new} & = & \mat{F}.
\end{eqnarray*}

\section{Convergence: continuous time systems}

We can split our considerations here to those concerning global, and those concerning local convergence properties. We start with global convergence.

\bigskip
\begin{mdframed}
\begin{theorem}[Global convergence]	\label{thm:globconv_cLQR}
Let us start the split optimal policy iteration with the initial guess $\mat{F}^0$ and with updating the $i^{\rm th}$ subsystem first, where $i$ is arbitrary. Then, if the pair $(\mat{A}^{(i)},\mat{B}^{(i)})$ is controllable, then
\begin{enumerate}[(i)]
	\item all the iterates $\mat{F}^k$, $k\ge 1$, stabilize the global system;
	\item the value functions $V^k$ corresponding to the $\mat{F}^k$ satisfy $V^{k+1}\le V^k$ pointwise; and
	\item $\mat{F}^k\to\mat{F}^{\rm opt}$ as $k\to\infty$, where $\mat{F}^{\rm opt}$ is the optimal feedback matrix of the original problem.
\end{enumerate}
\end{theorem}
\end{mdframed}

\bigskip

In order to prove this theorem, we break it down into several lemmas. Without mentioning it again, the controllability condition in the theorem is assumed to be valid throughout this section. First, we show the monotony of the value functions. Note that, since all the subproblems solved during the iteration are LQR problems themselves, $V^k(x) = x^{\rm T}\mat{P}^k x$ for some symmetric positive definite $\mat{P}^k$.
\begin{lemma}	\label{lem:monotony_cLQR}
We have $\mat{P}^{k+1} \preceq \mat{P}^k$ for $k\ge 0$, and $\mat{F}^k$ is stabilizing for every $k\ge 1$. Here, $\mat{A}\preceq\mat{B}$ for two symmetric matrices if and only if $\mat{B}-\mat{A}$ is symmetric positive semidefinite.
\end{lemma}
\begin{myproof}
Even though the iteration is based on the successive solution of subproblems, they are constructed such that the value function they optimize is the same object in each step: it is $\int_0^{\infty} x(t)^{\rm T}\mat{Q}x(t) + u(t)^{\rm T}\mat{R}u(t)\,dt$ with $u(t) = \mat{F}^kx(t)$ for the $k^{\rm th}$ iterate. By not changing the $\mat{F}^k$ in a step we would have $\mat{P}^{k+1} = \mat{P}^k$. But since in each step the optimal control (for the corresponding subproblem) is taken, we obtain $\mat{P}^{k+1} \preceq \mat{P}^k$ by optimality, and all corresponding feedbacks are stabilizing. The feedback matrix $\mat{F}^1$ is stabilizing by assumption.
\end{myproof}

Now, that we have a monotonically decreasing sequence of value functions, their convergence has to be established.
\begin{lemma}
Let $\{\mat{P}^k\}_{k\ge 1}$ a sequence of symmetric positive definite matrices, such that $\mat{P}^{k+1}\preceq \mat{P}^k$ for every $k\ge 1$. Then $\mat{P}^k\to\bar{\mat{P}}$ as $k\to\infty$ for some symmetric positive semidefinite $\bar{\mat{P}}$.
\end{lemma}
\begin{myproof}
(This elegant proof is due to Benedict Dingfelder.)
By assumption, $\{v^{\rm T}\mat{P}^k v\}_{k\ge 1}$ is a monotonically decreasing sequence of positive real numbers for any $v\in\R^m$, hence it converges, say, to $p_v\ge 0$. Setting $v=e_i$, the $i^{\rm th}$ canonical vector, we obtain $\mat{P}^k_{ii} \to p_{e_i}$ as $k\to\infty$. With $v_{ij} := e_i+e_j$ we have $v_{ij}^{\rm T}\mat{P}^kv_{ij} = \mat{P}^k_{ii} + \mat{P}^k_{jj} + 2\mat{P}^k_{ij}$, since all $\mat{P}^k$ are symmetric. The left hand side of this equation converges, just as the diagonal entries of the matrix, shown above. It follows that $\mat{P}_{ij}$ converges for every $i,j\in\{1,\ldots,m\}$, hence $\mat{P}^k$ converges element wise (and thus in any norm) to a matrix~$\bar{\mat{P}}$. Since the set of positive semidefinite matrices is a closed subspace of~$\R^{m\times m}$, the matrix~$\bar{\mat{P}}$ is contained in this subspace, as the limit of symmetric positive definite matrices.
\end{myproof}

Next is to show that the feedback matrices converge to a fixed point of the split optimal policy iteration.
\begin{lemma}	\label{lem:FP_cLQR}
It holds $\mat{F}^k\to\bar{\mat{F}}$ as $k\to\infty$, where $\bar{\mat{F}} = -\mat{R}^{-1}\mat{B}^{\rm T}\bar{\mat{P}}$ is a fixed point of the split optimal policy iteration.
\end{lemma}
\begin{myproof}
From the update assignment for the feedback matrices we have \linebreak[4] $\mat{\Pi}_i^{\rm T}\mat{F}^{k_i} = -{\mat{R}^{(i)}}^{-1}{\mat{B}^{(i)}}^{\rm T}\mat{P}^{k_i}$, where the $k_i$ are indices of the iteration steps where the $i^{\rm th}$ subsystem is updated. Thus, $\mat{\Pi}_i^{\rm T}\mat{F}^k\to -{\mat{R}^{(i)}}^{-1}{\mat{B}^{(i)}}^{\rm T}\bar{\mat{P}}$ as $k\to\infty$, because $\mat{P}^{k_i}\to\bar{\mat{P}}$ as $k_i\to\infty$, independently from~$i$. Since $\sum_i\mat{\Pi}_i\mat{\Pi}_i^{\rm T} = \mat{I}\in\R^{m\times m}$, we have
\begin{eqnarray*}
\bar{\mat{F}} & = & \sum_i\mat{\Pi}_i\mat{\Pi}_i^{\rm T}\bar{\mat{F}} \\
		   & = & -\sum_i\mat{\Pi}_i (\mat{\Pi}_i^{\rm T}\mat{R}\mat{\Pi}_i)^{-1}\mat{\Pi}_i^{\rm T}\mat{B}^{\rm T}\bar{\mat{P}} \\
		   & = & -\mat{R}^{-1}\mat{B}^{\rm T}\bar{\mat{P}}
\end{eqnarray*}
where the last equality follows from $\mat{R}$ being block diagonal.

In order to see that $\bar{\mat{F}}$ is a fixed point of the split optimal policy iteration, we investigate how $\mat{F}^{k+1}$ is obtained from~$\mat{F}^k$. The mapping $g^i:\mat{F}^k \mapsto \mat{F}^{k+1}$ can be decomposed into the following steps:
\begin{equation}
g^i: \mat{F}^k \mapsto (\mat{A}^{(i)},\mat{B}^{(i)},\mat{Q}^{(i)},\mat{R}^{(i)}) \mapsto \mat{P}^{k+1} \mapsto \mat{F}^{k+1}.
\label{eq:iter_map}
\end{equation}
The first and last mapping are obviously continuous (and arbitrarily often differentiable) by~\eqref{eq:matrixredef} and~\eqref{eq:contOptF}; the second one is even analytic~\cite{Del84}. Thus, $g^i$ is continuous  for every~$i$, and the equation $\mat{F}^{k+1} = g^i(\mat{F}^k)$ yields $\bar{\mat{F}} = g^i(\bar{\mat{F}})$ as $k\to\infty$.
\end{myproof}

We are now ready to show the theorem. It remains to be shown that the fixed point $\bar{\mat{F}}$ of the split optimal policy iteration coincides with the optimal solution of the initial problem, i.e.\ with~$\mat{F}^{\rm opt}$.

\begin{myproof}[Proof of Theorem~\ref{thm:globconv_cLQR}.]
Statements (i) and (ii) are shown in Lemma~\ref{lem:monotony_cLQR}. Statement~(iii) follows from Lemma~\ref{lem:FP_cLQR}, if we show that any fixed point of the split optimal policy iteration (more precisely the associated matrix~$\bar{\mat{P}}$) solves the continuous algebraic Riccati equation~\eqref{eq:CARE}, and by uniqueness of its solutions~\cite{Sont98} coincides with the optimal solution.

Since $\bar{\mat{P}}$ and $\bar{\mat{F}}$ are fixed points of the iteration, they solve~\eqref{eq:CARE} with the modified matrices~\eqref{eq:matrixredef} for every~$i$. This reads as
\begin{eqnarray*}
\mat{0} & = & \bar{\mat{P}}\mat{B}\mat{\Pi}_i\left(\mat{\Pi}_i^{\rm T}\mat{R}\mat{\Pi}_i\right)^{-1}\mat{\Pi}_i^{\rm T}\mat{B}^{\rm T}\bar{\mat{P}} - \bar{\mat{P}}\left(\mat{A}+\mat{B}\left(\mat{I}-\mat{\Pi}_i\mat{\Pi}_i^{\rm T}\right)\bar{\mat{F}}\right) \\
  &   & - \left(\mat{A}+\mat{B}\left(\mat{I}-\mat{\Pi}_i\mat{\Pi}_i^{\rm T}\right)\bar{\mat{F}}\right)^{\rm T}\bar{\mat{P}} - \mat{Q} - \bar{\mat{F}}^{\rm T}\left(\mat{I}-\mat{\Pi}_i\mat{\Pi}_i^{\rm T}\right)\mat{R}\left(\mat{I}-\mat{\Pi}_i\mat{\Pi}_i^{\rm T}\right)\bar{\mat{F}}.
\end{eqnarray*}
It is easy to see be the block diagonal structure of~$\mat{R}$ that
\[
\mat{\Pi}_i\left(\mat{\Pi}_i^{\rm T}\mat{R}\mat{\Pi}_i\right)^{-1}\mat{\Pi}_i^{\rm T} = \mat{\Pi}_i\mat{\Pi}_i^{\rm T}\mat{R}^{-1},
\]
and
\[
\left(\mat{I}-\mat{\Pi}_i\mat{\Pi}_i^{\rm T}\right) \mat{R} \left(\mat{I}-\mat{\Pi}_i\mat{\Pi}_i^{\rm T}\right) = \mat{R} \left(\mat{I}-\mat{\Pi}_i\mat{\Pi}_i^{\rm T}\right) = \left(\mat{I}-\mat{\Pi}_i\mat{\Pi}_i^{\rm T}\right) \mat{R}.
\]
Using these identities and $\bar{\mat{F}} = -\mat{R}^{-1}\mat{B}^{\rm T}\bar{\mat{P}}$ from Lemma~\ref{lem:FP_cLQR}, the above equation becomes
\begin{eqnarray*}
\mat{0} & = & \bar{\mat{P}}\mat{B}\mat{\Pi}_i\mat{\Pi}_i^{\rm T}\mat{R}^{-1}\mat{B}^{\rm T}\bar{\mat{P}} -
\bar{\mat{P}}\mat{A} - \mat{A}^{\rm T}\bar{\mat{P}} \\
  &   & + \bar{\mat{P}}\mat{B}\left(\mat{I}-\mat{\Pi}_i\mat{\Pi}_i^{\rm T}\right) \mat{R}^{-1} \mat{B}^{\rm T}\bar{\mat{P}} + \bar{\mat{P}}\mat{B}\mat{R}^{-1} \left(\mat{I}-\mat{\Pi}_i\mat{\Pi}_i^{\rm T}\right)\mat{B}^{\rm T}\bar{\mat{P}} \\
  &   & - \mat{Q} - \bar{\mat{P}}\mat{B}\left(\mat{I}-\mat{\Pi}_i\mat{\Pi}_i^{\rm T}\right)\mat{R}^{-1}\mat{B}^{\rm T}\bar{\mat{P}}
\end{eqnarray*}
By collecting all terms starting with $\bar{\mat{P}}\mat{B}$ we obtain
\begin{eqnarray*}
\mat{0} & = & \bar{\mat{P}}\mat{B}\left(\mat{\Pi}_i\mat{\Pi}_i^{\rm T}\mat{R}^{-1} +
\left(\mat{I}-\mat{\Pi}_i\mat{\Pi}_i^{\rm T}\right) \mat{R}^{-1} + \mat{R}^{-1} \left(\mat{I}-\mat{\Pi}_i\mat{\Pi}_i^{\rm T}\right) - \left(\mat{I}-\mat{\Pi}_i\mat{\Pi}_i^{\rm T}\right)\mat{R}^{-1} \right)\mat{B}^{\rm T}\bar{\mat{P}} \\
  &   & - \bar{\mat{P}}\mat{A} - \mat{A}^{\rm T}\bar{\mat{P}} - \mat{Q} \\
	& = & \bar{\mat{P}}\mat{B}\mat{R}^{-1}\mat{B}^{\rm T}\bar{\mat{P}} - \bar{\mat{P}}\mat{A} - \mat{A}^{\rm T}\bar{\mat{P}} - \mat{Q},
\end{eqnarray*}
hence $\bar{\mat{P}}$ solves~\eqref{eq:CARE}, and this proves the claim.
\end{myproof}

\vspace*{-1em}
\begin{remark}
It is enough to have controllability with respect to one subsystem: if there is a non-controllable subsystem, we just jump to the next without changing the feedback matrix as long as we find a subsystem with respect to which one has controllability.
\end{remark}

\smallskip
Locally, around the optimal solution, the iteration can be shown to converge quickly:

\pagebreak[4]
\begin{mdframed}
\begin{theorem}[Local speed of convergence]	\label{thm:convspeed_cLQR}
The split optimal policy iteration for continuous time LQR problem converges locally quadratically; i.e.\ for $\left\|\mat{F}^k-\mat{F}^{\rm opt}\right\|$ sufficiently small we have
\[
\left\| \mat{F}^{k+n} - \mat{F}^{\rm opt}\right\| \le C \left\|\mat{F}^k - \mat{F}^{\rm opt}\right\|^2
\]
for some $C>0$ which is independent of the $\mat{F}^k$.
\end{theorem}
\end{mdframed}

\bigskip

\begin{myproof}
The iteration is a consecutive application of the mappings $g^i$ form~\eqref{eq:iter_map} for $i$ cyclically sweeping over the subsystem indices $\{1,\ldots,n\}$. One cycle can be written as $g:=g^{j-1}\circ g^{j-2}\circ\cdots\circ g^{j+1}\circ g^j$ for some starting subsystem index~$j$. The global convergence of the successive iterates of~$g$ to $\mat{F}^{\rm opt}$ from~\eqref{eq:contOptF} was established above. Since $g(\mat{F}^{\rm opt}) = \mat{F}^{\rm opt}$, and $g$ is arbitrary often differentiable, the local convergence depends on the spectral radius of the derivative of~$g$ at~$\mat{F}^{\rm opt}$.

We will show local quadratic convergence in the sense that
\[
\left\| g(\mat{F}) - \mat{F}^{\rm opt}\right\| \le C \left\|\mat{F} - \mat{F}^{\rm opt}\right\|^2
\]
for any $\mat{F}$ from a given neighborhood of the optimal solution, and some $C>0$ independent of~$\mat{F}$. Since $g$ is a composition of the $g^i$, which all update different rows of $\mat{F}$ during the iteration, it suffices to show quadratic convergence for the rows updated by a specific $g^i$ for an arbitrary~$i$, and the claim follows. Hence, fix $i\in\{1\,\ldots,n\}$ and consider the map~$g^i$ from~\eqref{eq:iter_map}. Quadratic convergence follows if $Dg^i(\mat{F}^{\rm opt}) =0$. Since the map $\mat{P}^{k+1}\mapsto \mat{F}^{k+1}$ is independent of $\mat{F}^k$, we have $Dg^i(\mat{F}^{\rm opt})=0$ if $\partial\mat{P}^{k+1}/\partial\mat{F}^k\,(\mat{F}^{\rm opt})=0$.

To see this, note that $\mat{P}^{k+1}$ is the unique symmetric positive definite solution \linebreak[4] of~$r^i(\mat{P}^{k+1},\mat{F}^k)=\mat{0}$, where
\[
r^i(\mat{P},\mat{F}) := \mat{P}\mat{B}^{(i)}{\mat{R}^{(i)}}^{-1}{\mat{B}^{(i)}}^{\rm T}\mat{P} - \mat{P}\mat{A}^{(i)} - {\mat{A}^{(i)}}^{\rm T}\mat{P} - \mat{Q}^{(i)},
\]
with the definitions from~\eqref{eq:matrixredef}, and the matrices $\mat{A}^{(i)}$ and $\mat{Q}^{(i)}$ depending on~$\mat{F}$. By the implicit function theorem one has for $\mat{P}(\mat{F})$ satisfying $r^i(\mat{P}(\mat{F}),\mat{F})=\mat{0}$ that
\[
\frac{\partial\mat{P}}{\partial\mat{F}}(\mat{F}) = -\left(\frac{\partial r^i}{\partial \mat{P}}(\mat{P},\mat{F})\right)^{-1}\,\frac{\partial r^i}{\partial\mat{F}}(\mat{P},\mat{F}),
\]
in some neighborhood of any point such that $\partial r^i/\partial\mat{P}$ is invertible. We will show \linebreak[4] ${\tfrac{\partial r^i}{\partial \mat{F}}(\mat{P}^{\rm opt},\mat{F}^{\rm opt})=\mat{0}}$, implying the claim.

Differentiating $r^i$ yields
\begin{equation}
\begin{aligned}
\frac{\partial r^i}{\partial \mat{F}}(\mat{P},\mat{F})\cdot \mat{\Delta} &= -\mat{P} \mat{B} \left(\mat{I}-\mat{\Pi}_i\mat{\Pi}_i^{\rm T}\right)\mat{\Delta} - \mat{F}^{\rm T}\left(\mat{I}-\mat{\Pi}_i\mat{\Pi}_i^{\rm T}\right)\mat{R} \left(\mat{I}-\mat{\Pi}_i\mat{\Pi}_i^{\rm T}\right) \mat{\Delta} \\
&{\phantom{=}} \ - \mat{\Delta}^{\rm T} \left(\mat{I}-\mat{\Pi}_i\mat{\Pi}_i^{\rm T}\right) \mat{B}^{\rm T}\mat{P} - \mat{\Delta}^{\rm T}\left(\mat{I}-\mat{\Pi}_i\mat{\Pi}_i^{\rm T}\right) \mat{R} \left(\mat{I}-\mat{\Pi}_i\mat{\Pi}_i^{\rm T}\right) \mat{F},
\end{aligned}
\label{eq:diffCARE}
\end{equation}
for every $\mat{\Delta}\in\R^{r\times m}$. Substituting $\mat{F}^{\rm opt} = -\mat{R}^{-1}\mat{B}^{\rm T}\mat{P}^{\rm opt}$, we get
\begin{eqnarray*}
{\mat{F}^{\rm opt}}^{\rm T}\left(\mat{I}-\mat{\Pi}_i\mat{\Pi}_i^{\rm T}\right)\mat{R} \left(\mat{I}-\mat{\Pi}_i\mat{\Pi}_i^{\rm T}\right) & = & -\mat{P}^{\rm opt}\mat{B}\mat{R}^{-1}\left(\mat{R} - \mat{\Pi}_i\mat{\Pi}_i^{\rm T}\mat{R}\right)\left(\mat{I}-\mat{\Pi}_i\mat{\Pi}_i^{\rm T}\right) \\
& = & -\mat{P}^{\rm opt}\mat{B}\left(\mat{I}-\mat{\Pi}_i\mat{\Pi}_i^{\rm T}\right)^2 \\
& = & -\mat{P}^{\rm opt}\mat{B}\left(\mat{I}-\mat{\Pi}_i\mat{\Pi}_i^{\rm T}\right),
\end{eqnarray*}
since $\mat{R}^{-1}\mat{\Pi}_i\mat{\Pi}_i^{\rm T}\mat{R} = \mat{\Pi}_i\mat{\Pi}_i^{\rm T}$ by the block diagonal structure of~$\mat{R}$, and $(\mat{\Pi}_i\mat{\Pi}_i^{\rm T})^2 = \mat{\Pi}_i\mat{\Pi}_i^{\rm T}$. Using this in~\eqref{eq:diffCARE} shows ${\tfrac{\partial r^i}{\partial \mat{F}}(\mat{P}^{\rm opt},\mat{F}^{\rm opt})=\mat{0}}$. This concludes the proof.
\end{myproof}

\begin{remark}
The proofs show that the only assumption we need on the structure of the matrices $\mat{A}, \mat{B}, \mat{Q}$ and $\mat{R}$ is that $\mat{R}$ is block diagonal. There are no restrictions on the state, input and state cost matrices.
\end{remark}

\section{Convergence: discrete time systems}

\def\hI{\hat{\mat{I}}_i}

We can proceed in analogous manner as we did in the previous section. Again, we start with global convergence, where exactly the same claim holds.

\bigskip

\begin{mdframed}
\begin{theorem}[Global convergence]	\label{thm:globconv_dLQR}
Let us start the split optimal policy iteration with the initial guess $\mat{F}^0$ and with updating the $i^{\rm th}$ subsystem first, where $i$ is arbitrary. Then, if the pair $(\mat{A}^{(i)},\mat{B}^{(i)})$ is controllable, then
\begin{enumerate}[(i)]
	\item all the iterates $\mat{F}^k$, $k\ge 1$, stabilize the global system;
	\item the value functions $V^k$ corresponding to the $\mat{F}^k$ satisfy $V^{k+1}\le V^k$ pointwise; and
	\item $\mat{F}^k\to\mat{F}^{\rm opt}$ as $k\to\infty$, where $\mat{F}^{\rm opt}$ is the optimal feedback matrix of the original problem.
\end{enumerate}
\end{theorem}
\end{mdframed}

\bigskip

Claims (i) and (ii) of this theorem are proven the same way as for Theorem~\ref{thm:globconv_cLQR}, the adaptation of lemmas \ref{lem:monotony_cLQR}--\ref{lem:FP_cLQR} is done by just replacing the formulas for the continuous time system with the ones for the discrete time system. However, claim (iii) does need a little bit more work, and in order to improve the readability of the proof below, we prove the most technical part in the following lemma---which can be skipped if the reader is not interested in the details.

\begin{lemma}	\label{lem:tech_mat_id}
Let $\mat{S}$ and $\mat{R}$ be symmetric positive definite matrices, $\mat{R}$ being also block diagonal. Set $\hI = \mat{I}-\mat{
\Pi}_i\mat{\Pi}_i^{\rm T}$, $\mat{Z} = \left(\mat{R}+\mat{S}\right)^{-1}$, $\mat{Z}_i = \mat{\Pi}_i\left(\mat{\Pi}_i^{\rm T}\left(\mat{R}+\mat{S}\right)\mat{\Pi}_i\right)^{-1}\mat{\Pi}_i^{\rm T}$. Then the following equations hold:
\begin{equation}
\mat{Z}^{-1}\mat{Z}_i\mat{Z}^{-1} - \hI\mat{S}\mat{Z}_i\mat{Z}^{-1} - \mat{Z}^{-1}\mat{Z}_i\mat{S}\hI + \hI\mat{S}\mat{Z}_i\mat{S}\hI = \mat{\Pi}_i\mat{\Pi}_i^{\rm T}\left(\mat{R}+\mat{S}\right)\mat{\Pi}_i\mat{\Pi}_i^{\rm T}
\label{eq:tech_mat_id1}
\end{equation}
\begin{equation}
\mat{Z}^{-1} = \hI\mat{Z}^{-1} + \mat{Z}^{-1}\hI - \hI\left(\mat{R}+\mat{S}\right)\hI + \mat{\Pi}_i\mat{\Pi}_i^{\rm T}\left(\mat{R}+\mat{S}\right)\mat{\Pi}_i\mat{\Pi}_i^{\rm T}
\label{eq:tech_mat_id2}
\end{equation}
\end{lemma}

\begin{myproof}
To prove \eqref{eq:tech_mat_id1}, consider the second and fourth terms on the left hand side:
\begin{eqnarray*}
\hI\mat{S}\mat{Z}_i\mat{S}\hI - \hI\mat{S}\mat{Z}_i\mat{Z}^{-1} & = & \hI\mat{S}\mat{Z}_i\left(\mat{S}\hI-\mat{R}-\mat{S}\right) \\
	& = & -\hI\mat{S}\mat{Z}_i\left(\mat{R}+\mat{S}\right)\mat{\Pi}_i\mat{\Pi}_i^{\rm T}
\end{eqnarray*}
where we used $\mat{Z}_i\mat{R} = \mat{Z}_i\mat{R}\mat{\Pi}_i\mat{\Pi}_i^{\rm T}$, because~$\mat{Z}_i\mat{R}$ has non-zero blocks only on its~$i^{\rm th}$ block row. Analogously, the first and third terms yield:
\[
\mat{Z}^{-1}\mat{Z}_i\mat{Z}^{-1} - \mat{Z}^{-1}\mat{Z}_i\mat{S}\hI = \left(\mat{R}+\mat{S}\right) \mat{Z}_i \left(\mat{R}+\mat{S}\right)\mat{\Pi}_i\mat{\Pi}_i^{\rm T}.
\]
Adding up these two equations gives us that the left hand side of \eqref{eq:tech_mat_id1} is equal to
\[
\left(\mat{R}+\mat{S} - \hI\mat{S}\right) \mat{Z}_i \left(\mat{R}+\mat{S}\right)\mat{\Pi}_i\mat{\Pi}_i^{\rm T} = \mat{\Pi}_i\mat{\Pi}_i^{\rm T}\left(\mat{R}+\mat{S}\right) \mat{Z}_i \left(\mat{R}+\mat{S}\right)\mat{\Pi}_i\mat{\Pi}_i^{\rm T}.
\]
Now we note that $\mat{\Pi}_i\mat{\Pi}_i^{\rm T}\left(\mat{R}+\mat{S}\right)$ coincides with $\mat{R}+\mat{S}$ in its $i^{\rm th}$ block row and is zero elsewhere, and $\left(\mat{R}+\mat{S}\right)\mat{\Pi}_i\mat{\Pi}_i^{\rm T}$ coincides with $\mat{R}+\mat{S}$ in its $i^{\rm th}$ block column and is zero elsewhere. The matrix $\mat{Z}_i$ is the inverse of the $(i,i)^{\rm th}$ block of $\mat{R}+\mat{S}$ in its $(i,i)^{\rm th}$ block and is zero elsewhere. Hence the product of these matrices is exactly the right hand side of \eqref{eq:tech_mat_id1}, which proves the equality.

To see that  \eqref{eq:tech_mat_id2} holds, note that
\begin{itemize}
\item $\hI\left(\mat{R}+\mat{S}\right)$ is the matrix $\left(\mat{R}+\mat{S}\right)$ with its $i^{\rm th}$ block row set to zero,
\item $\left(\mat{R}+\mat{S}\right)\hI$ is the matrix $\left(\mat{R}+\mat{S}\right)$ with its $i^{\rm th}$ block column set to zero,
\item $\hI\left(\mat{R}+\mat{S}\right)\hI$ is the matrix $\left(\mat{R}+\mat{S}\right)$ with its $i^{\rm th}$ row and column set to zero, and
\item $\mat{\Pi}_i\mat{\Pi}_i^{\rm T}\left(\mat{R}+\mat{S}\right)\mat{\Pi}_i\mat{\Pi}_i^{\rm T}$ is the matrix $\left(\mat{R}+\mat{S}\right)$ with everything apart from its $(i,i)^{\rm th}$ block set to zero.
\end{itemize}
Combining the terms on the right hand side as they are gives exactly the matrix $\mat{R}+\mat{S}$, i.e.\ $\mat{Z}^{-1}$.
\end{myproof}

Now everything is set up to prove Theorem~\ref{thm:globconv_dLQR}.
\begin{myproof}[Proof of Theorem~\ref{thm:globconv_dLQR}.]
Statements (i) and (ii) are shown analogously as before. Statement~(iii) follows from Lemma~\ref{lem:FP_cLQR} (more precisely, from its discrete counterpart), if we show that any fixed point of the split optimal policy iteration (more precisely the associated matrix~$\bar{\mat{P}}$) solves the discrete algebraic Riccati equation~\eqref{eq:DARE}, and by uniqueness of its solutions~\cite{Sont98} coincides with the optimal solution.

Since $\bar{\mat{P}}$ and $\bar{\mat{F}}$ are fixed points of the iteration, they solve~\eqref{eq:DARE} with the modified matrices~\eqref{eq:matrixredef} for every~$i$. This reads as, by using the notation of the previous lemma:
\begin{eqnarray*}
\bar{\mat{P}} & = & \left(\mat{A}+\mat{B}\hI\bar{\mat{F}}\right)^{\rm T}\bar{\mat{P}}\left(\mat{A}+\mat{B}\hI\bar{\mat{F}}\right) - \left(\mat{A}+\mat{B}\hI\bar{\mat{F}}\right)^{\rm T}\bar{\mat{P}}\mat{B}\mat{Z}_i\mat{B}^{\rm T}\bar{\mat{P}}\left(\mat{A}+\mat{B}\hI\bar{\mat{F}}\right) + \mat{Q} + \bar{\mat{F}}^{\rm T}\hI\mat{R}\hI\bar{\mat{F}} \\
	& = & \mat{A}^{\rm T}\bar{\mat{P}}\mat{A} - \mat{A}^{\rm T}\bar{\mat{P}}\mat{B}\mat{M}\mat{B}^{\rm T}\bar{\mat{P}}\mat{A} + \mat{Q},
\end{eqnarray*}
with
\[
\mat{M} := \mat{Z}\left( \hI\mat{Z}^{-1} + \mat{Z}^{-1}\hI - \hI\left(\mat{R}+\mat{S}\right)\hI + \mat{Z}^{-1}\mat{Z}_i\mat{Z}^{-1} - \hI\mat{S}\mat{Z}_i\mat{Z}^{-1} - \mat{Z}^{-1}\mat{Z}_i\mat{S}\hI + \hI\mat{S}\mat{Z}_i\mat{S}\hI \right)\mat{Z}\, ,
\]
where $\mat{S} = \mat{B}^{\rm T}\bar{\mat{P}}\mat{B}$. The above equation is identical with the full discrete algebraic Riccati equation~\eqref{eq:DARE} if $\mat{M} = \mat{Z}$. This identity is shown in Lemma~\ref{lem:tech_mat_id}. Hence we have shown that the fixed point $\bar{\mat{P}}$ of the split optimal policy iteration also solves~\eqref{eq:DARE}, and by uniqueness of the solutions it is identical with the solution of this full equation. This concludes the proof.
\end{myproof}

Having shown global convergence, we turn to the analysis of the convergence speed. Quite surprisingly, it turns out that for the discrete time case the convergence is not local quadratic in general, but is in a non-trivial relation with the coupling strength between the subsystems. The weaker the coupling between the subsystems, the faster the convergence of the split optimal policy iteration; and if no coupling is present, the convergence is locally quadratic.

Let us recall, that just as in the continuous times case, one iteration cycle can be written as $g:=g^{j-1}\circ g^{j-2}\circ\cdots\circ g^{j+1}\circ g^j$ (if we start the iteration at the $j^{\rm th}$ subsystem), where the $g^i$, $i=1,\ldots,n$, are the maps describing the update of the feedback matrix~$\mat{F}$. Since $\mat{B}^{(i)}$ and $\mat{R}^{(i)}$ do not depend on~$\mat{F}$, the map~$g^i$ can be decomposed into factors:
\[
g^i:\ \mat{F}\mapsto \big(\mat{A}^{(i)},\mat{Q}^{(i)}\big)\mapsto \big(\mat{P}^{\rm new},\mat{A}^{(i)}\big)\mapsto \mat{F}^{\rm new},
\]
where the first mapping encodes the modified system matrices due to \eqref{eq:matrixredef}, the second realizes the solution of the Riccati equation~\eqref{eq:DARE}, and the third performs the update of the feedback matrix as described in section~\ref{sec:update_general}.

\bigskip

\begin{mdframed}
\begin{theorem}[Local speed of convergence] 	\label{thm:convspeed_dLQR}
The split optimal policy iteration for discrete time LQR problem converges locally linearly, i.e.\ there is a $\varrho>0$ such that
\[
\left\| \mat{F}^{k+n} - \mat{F}^{\rm opt}\right\| \le \varrho \left\|\mat{F}^k - \mat{F}^{\rm opt}\right\|
\]
The convergence rate~$\varrho$ is governed by the norm of the iteration matrix
\[
Dg(\mat{F}^{\rm opt})=Dg^{j-1}(\mat{F}^{\rm opt})\cdot Dg^{j-2}(\mat{F}^{\rm opt})\cdots Dg^{j+1}(\mat{F}^{\rm opt})\cdot Dg^j(\mat{F}^{\rm opt}),
\]
where
\begin{equation}
Dg^i(\mat{F}^{\rm opt}) = -\mat{\Pi}_i \left(\mat{R}^{(i)} + {\mat{B}^{(i)}}^{\rm T} \mat{P}^{\rm opt}\mat{B}^{(i)}\right)^{-1}{\mat{B}^{(i)}}^{\rm T}\mat{P}^{\rm opt}\mat{B}\hI\,.
\label{eq:sens_mat}
\end{equation}
\end{theorem}
\end{mdframed}

\bigskip

\begin{myproof}
The strategy of the proof follows that of Theorem~\ref{thm:convspeed_cLQR}, although it will turn out that $Dg^i(\mat{F}^{\rm opt})\neq \mat{0}$ in general. Using the above factorization of $g^i$ we can give an expression for its derivative. We use the abbreviation $\mat{Z}_i = \mat{\Pi}_i\left(\mat{R}^{(i)} + {\mat{B}^{(i)}}^{\rm T} \mat{P}^{\rm new}\mat{B}^{(i)}\right)^{-1}\mat{\Pi}_i^{\rm T}$.
\begin{eqnarray*}
Dg^i(\mat{F}) & = & \frac{\partial \mat{F}^{\rm new}}{\partial \mat{F}} \\
    & = & \frac{\partial \mat{F}^{\rm new}}{\partial \mat{P}^{\rm new}}\frac{\partial \mat{P}^{\rm new}}{\partial \mat{F}} + \frac{\partial \mat{F}^{\rm new}}{\partial \mat{A}^{(i)}}\frac{\partial \mat{A}^{(i)}}{\partial \mat{F}} \\
	& = & \mat{Z}_i \mat{B}^{\rm T}\frac{\partial \mat{P}^{\rm new}}{\partial \mat{F}} \mat{B} \mat{Z}_i \mat{B}^{\rm T}\mat{P}^{\rm new}\mat{A}^{(i)} - \mat{Z}_i \mat{B}^{\rm T}\frac{\partial \mat{P}^{\rm new}}{\partial \mat{F}}\mat{A}^{(i)} - \mat{Z}_i \mat{B}^{\rm T}\mat{P}^{\rm new}\frac{\partial \mat{A}^{(i)}}{\partial \mat{F}} \\
	& = & \mat{Z}_i \mat{B}^{\rm T}\left[ \frac{\partial \mat{P}^{\rm new}}{\partial \mat{F}}\mat{B}\mat{Z}_i \mat{B}^{\rm T}\mat{P}^{\rm new}\mat{A}^{(i)} - \frac{\partial \mat{P}^{\rm new}}{\partial \mat{F}}\mat{A}^{(i)} - \mat{P}^{\rm new}\frac{\partial \mat{A}^{(i)}}{\partial \mat{F}} \right]
\end{eqnarray*}
where the third equality follows directly from the update formula (section~\ref{sec:update_general}) and the chain rule.

Next, we are going to determine $\frac{\partial \mat{P}^{\rm new}}{\partial \mat{F}}$ via the implicit function theorem. Note that $\mat{P}^{\rm new}$ is defined by $r^i(\mat{P}^{\rm new},\mat{F})=\mat{0}$, where
\[
r^i(\mat{P},\mat{F}):= {\mat{A}^{(i)}}^{\rm T}\mat{P}\mat{A}^{(i)} - {\mat{A}^{(i)}}^{\rm T}\mat{P}\mat{B}^{(i)}\left(\mat{R}^{(i)}+{\mat{B}^{(i)}}^{\rm T}\mat{P}\mat{B}^{(i)}\right)^{-1}{\mat{B}^{(i)}}^{\rm T}\mat{P}\mat{A}^{(i)} + \mat{Q}^{(i)} - \mat{P}\,.
\]
Let us recall that $\mat{P}^{\rm new}$ is analytic in every parameter of the Riccati equation~\cite{Del84}, hence there is no obstacle in our way to use the implicit function theorem. Just as in the proof of Theorem~\ref{thm:convspeed_cLQR} we will show that $\frac{\partial r^i}{\partial\mat{F}}(\mat{F}^{\rm opt})=\mat{0}$, implying that $\frac{\partial \mat{P}^{\rm new}}{\partial \mat{F}}(\mat{F}^{\rm opt})=\mat{0}$.

To this end, we substitute in $r^i$ the matrices from~\eqref{eq:matrixredef}, differentiate it with respect to~$\mat{F}$ at $\mat{F}^{\rm opt}$ from~\eqref{eq:discrOptF}. Then, we pull out the common factor and simplify the result by using the notation from Lemma~\ref{lem:tech_mat_id} and Theorem~\ref{thm:globconv_dLQR}. What we get is
\begin{equation}
\frac{\partial r^i}{\partial\mat{F}}(\mat{F}^{\rm opt})\cdot\mat{\Delta} = \mat{A}^{\rm T}\mat{P}^{\rm opt}\mat{B}\left[\hI - \mat{Z}\hI\mat{Z}^{-1}\hI - \mat{Z}_i\mat{S}\hI + \mat{Z}\hI\mat{S}\mat{Z}_i\mat{S}\hI\right]\mat{\Delta} +\circledast\,,
\label{eq:diffDARE}
\end{equation}
where $\circledast$ is the transposed of the first term. Consider the last two terms on the right hand side of this expression in the brackets~$[\ldots ]$. Noting that $\mat{Z}^{-1} = \mat{R}+\mat{S}$, and using the shorthand $\mat{S}_{ii}$ and $\mat{R}_{ii}$ for the $(i,i)^{\rm th}$ blocks of $\mat{S}$ and $\mat{R}$, respectively, we obtain
\begin{eqnarray*}
- \mat{Z}_i\mat{S}\hI + \mat{Z}\hI\mat{S}\mat{Z}_i\mat{S}\hI & = & \mat{Z}\left(\hI\mat{S}-\mat{R}-\mat{S}\right)\mat{Z}_i\mat{S}\hI \\
 & = & \mat{Z}\left(-\mat{\Pi}_i\mat{\Pi}_i^{\rm T}\mat{S}-\mat{R}\right)\mat{\Pi}_i(\mat{R}_{ii}+\mat{S}_{ii})^{-1}\mat{\Pi}_i^{\rm T}\mat{S}\hI \\
  & = & -\mat{Z}\mat{\Pi}_i(\mat{S}_{ii}+\mat{R}_{ii})(\mat{S}_{ii}+\mat{R}_{ii})^{-1}\mat{\Pi}_i^{\rm T}\mat{S}\hI \\
  & = & -\mat{Z}\mat{\Pi}_i\mat{\Pi}_i^{\rm T}\mat{S}\hI,
\end{eqnarray*}
where we used in the third equality that $\mat{R}\mat{\Pi}_i = \mat{\Pi}_i\mat{R}_{ii}$ due to the fact that $\mat{R}$ is block diagonal. With this, returning to the bracketed expression in~\eqref{eq:diffDARE}, we have
\begin{eqnarray*}
[\ldots ] & = & \left(\mat{I} - \mat{Z}\left(\hI\mat{Z}^{-1} + \mat{\Pi}_i\mat{\Pi}_i^{\rm T}\mat{S}\right)\right)\hI \\
 & = & \left(\mat{I} - \mat{Z}\left(\hI\mat{R}+\mat{S}\right)\right)\hI \\
 & = & \mat{Z}\left(\mat{R}+\mat{S} - \hI\mat{R}-\mat{S}\right)\hI \\
 & = & \mat{Z}\mat{\Pi}_i\mat{\Pi}_i^{\rm T}\mat{R}\hI \\
 & = & \mat{0},
\end{eqnarray*}
since $\mat{\Pi}_i\mat{\Pi}_i^{\rm T}\mat{R}\hI=\mat{0}$, again, due to the block diagonal structure of~$\mat{R}$.

Returning to the formula for $Dg^i$, only the last term in the brackets does not vanish. We know from Theorem~\ref{thm:globconv_dLQR} that if $\mat{F} = \mat{F}^{\rm opt}$, then $\mat{P}^{\rm new} = \mat{P}^{\rm opt}$, hence we obtain
\[
Dg^i(\mat{F}^{\rm opt}) = -\mat{Z}_i\mat{B}^{\rm T}\mat{P}^{\rm opt}\mat{B}\hI\,.
\]
\end{myproof}

\begin{remark}[Rate of convergence]
If we assume distributed actuation (i.e.\ that the control input of the $i^{\rm th}$ subsystem has only a direct influence on the $i^{\rm th}$ subsystem itself), the input matrix~$\mat{B}$ is block diagonal with blocks $\mat{B}_{ii}$. Assuming a similar block decomposition according to the subsystem structure of the other matrices involved, \eqref{eq:sens_mat} yields that all block rows of $Dg^i(\mat{F}^{\rm opt})$ are zero apart from the $i^{\rm th}$ one, and its $(i,j)^{\rm th}$ block is
\begin{equation}
\big(\mat{R}_{ii}+\mat{B}_{ii}^{\rm T}\mat{P}_{ii}^{\rm opt}\mat{B}_{ii}\big)^{-1}\mat{B}_{ii}^{\rm T}\mat{P}_{ij}^{\rm opt}\mat{B}_{jj}
\label{eq:sens_block}
\end{equation}
for $j\neq i$, and is the zero matrix if $j=i$.

This can be interpreted in terms of coupling. If the system matrix $\mat{A}$ is block diagonal as well, i.e.\ the whole system consists of individual non-coupled subsystems, then $\mat{P}^{\rm opt}$ is block diagonal too and the matrix~\eqref{eq:sens_block} is zero for every~$j$. If the system is weakly coupled, i.e.\ the off-diagonal blocks of $\mat{A}$ are much smaller than its diagonal blocks, perturbation arguments~\cite{KoPeCh86} show that in general this property carries over to $\mat{P}^{\rm opt}$ as well.\footnote{In~\cite{KoJu14} the important case is studied where weak coupling in~$\mat{A}$ does not carry over to weak coupling in~$\mat{P}^{\rm opt}$.} Thus, \eqref{eq:sens_block} suggests that for weakly coupled systems we should expect the split optimal policy iteration to converge rapidly because $\varrho\ll 1$.
\end{remark}

{
\small
\bibliography{References}
\bibliographystyle{alpha}
}
\end{document}